\documentclass{amsart}

\usepackage{caption}

\usepackage{amssymb}
\usepackage[]{latexsym}

\usepackage{amsthm}

\usepackage{color}
\usepackage[usenames,dvipsnames,svgnames,table]{xcolor}
\usepackage[all,dvips]{xy}
\usepackage[dvips]{graphicx}

\usepackage{multirow}

\newtheorem{thm}{Theorem}
\newtheorem{cor}[thm]{Corollary}
\newtheorem{lemma}[thm]{Lemma}

\newtheorem{prop}[thm]{Proposition}

\theoremstyle{remark}
\newtheorem{rem}[thm]{Remark}
\newtheorem{exmp}[thm]{Example}

\theoremstyle{definition}

\newenvironment{pf}{\par\noindent{\bf Proof.}\enspace\ignorespaces}{\qed\par\par}

\def\qed{\hfill $\Box$}

\newcommand{\jac}[1]{{\rm Jac}( {#1} ) }
\newcommand{\Jac}{\operatorname{Jac}}

\newcommand{\pmodd}[1]{\,({\rm mod\,}{#1}) }

\newcommand{\im}{\mbox{Im}}

\newcommand{\Sel}{\mbox{Sel}}

\newcommand{\bQ}{{\mathbb{Q}}}

\newcommand{\bP}{{\mathbb{P}}}
\newcommand{\bZ}{{\mathbb{Z}}}
\newcommand{\bN}{{\mathbb{N}}}
\newcommand{\bA}{{\mathbb{A}}}

\newcommand{\cO}{{\mathcal{O}}}
\newcommand{\cZ}{{\mathcal{Z}}}
\newcommand{\cS}{{\mathcal{S}}}
\newcommand{\cT}{{\mathcal{T}}}
\newcommand{\cF}{{\mathcal{F}}}

\begin{document}

\title[On a conjecture of Rudin on squares in Arithmetic Progressions]{On a conjecture of Rudin\\ on squares in Arithmetic Progressions}

\author[Enrique Gonz\'alez-Jim\'enez]{Enrique Gonz\'alez-Jim\'enez}
\address{Universidad Aut{\'o}noma de Madrid, Departamento de Matem{\'a}ticas and Instituto de Ciencias Matem{\'a}ticas (ICMat), Madrid, Spain}
\email{enrique.gonzalez.jimenez@uam.es}
\urladdr{http://www.uam.es/enrique.gonzalez.jimenez}
\author[Xavier Xarles]{Xavier Xarles}
\address{Departament de Matem\`atiques\\Universitat Aut\`onoma de
Barcelona\\08193 Bellaterra, Barcelona, Spain}
\email{xarles@mat.uab.cat}

\subjclass[2010]{Primary: 11G30, 11B25,11D45; Secondary: 14H25}
\keywords{Arithmetic progressions, squares, covering collections, elliptic curve Chabauty}
\thanks{The first author was partially supported by the grant MTM2009--07291.
The second author was partially supported by the grant
MTM2006--11391.}


\begin{abstract}
Let $Q(N;q,a)$ denotes the number of squares in the arithmetic
progression $qn+a$, for $n=0, 1,\dots,N-1$, and let $Q(N)$ be the
maximum of $Q(N;q,a)$ over all non-trivial arithmetic progressions
$qn + a$.  Rudin's conjecture asserts that $Q(N)=O(\sqrt{N})$, and
in its stronger form that $Q(N)=Q(N;24,1)$ if $N\ge 6$. We prove
the conjecture above for $6\le N\le 52$. We even prove that the
arithmetic progression $24n+1$ is the only one, up to equivalence,
that contains $Q(N)$ squares for the values of $N$ such that
$Q(N)$ increases, for $7\le N\le 52$ (hence, for
$N=8,13,16,23,27,36,41$ and $52$). This allow us to assert, what we
have called Super--Strong Rudin's Conjecture: let be
$N=\mathcal{GP}_k+1\ge 8$ for some integer $k$, where
$\mathcal{GP}_k$ is the $k$-th generalized pentagonal number, then
$Q(N)=Q(N;q,a)$ with $\gcd(q,a)$ squarefree and $q> 0$ if and only if $(q,a)=(24,1)$.
\end{abstract}

\maketitle


\section{Introduction}

A well--known result by Fermat states that no four squares
in arithmetic progression over $\bZ$ exist. This result may be
reformulated in the following form: in four consecutive terms of a non-constant arithmetic progression, there are at most three squares. Hence, it is natural to ask how many squares there may be in $N$
consecutive terms of a non-constant arithmetic progression.

Following Bombieri, Granville and Pintz \cite{Bombieri-Granville-Pintz1992}, given $q$ and
$a$ integers, $q\ne 0$, we denote by $Q(N;q,a)$ the number of
squares in the arithmetic progression $qn+a$, for $n=0,1,\dots,N-1$
(there is a slight difference between our notation and the one in
\cite{Bombieri-Granville-Pintz1992}, since our arithmetic progressions begin with $i=0$ instead
of $i=1$). Denote by $Q(N)$ the maximum of $Q(N;q,a)$ over
all non-trivial arithmetic progressions. Notice that Fermat's
result is equivalent to $Q(4)=3$.

As a consequence of Fermat's result and of its own result on
arithmetic progressions, Szemer{\'e}di \cite{Szemeredi1975} proved an old
Erd\"os \cite{Erdos1963} conjecture: $Q(N)=o(N)$. This bound was
improved by Bombieri, Granville and Pintz \cite{Bombieri-Granville-Pintz1992} to
$Q(N)=O(N^{2/3+o(1)})$, and by Bombieri and Zannier \cite{Bombieri-Zannier2002} to
$Q(N)=O(N^{3/5+o(1)})$. Moreover, the so called Rudin's conjecture (\cite[end of \S 4.6]{Rudin1960})
asserts that $Q(N)=O(\sqrt{N})$, and in its stronger form that
(what we called Strong Rudin's Conjecture):
$$
Q(N)=Q(N;24,1)=\sqrt{\frac{8}{3} N}+O(1)\qquad \mbox{if $N\ge 6$.}
$$
Notice that $Q(5;24,1)=3$, but $Q(5;120,49)=4$ (since $7^2=49, 13^2=169, 17^2=289, 409, 23^2=529$)\footnote{It has been
proved in \cite{Gonzalez-Jimenez-Xarles2013} that the first quadratic number field where
there are five squares in arithmetic progression is
$\bQ(\sqrt{409})$; and  that the unique non-constant arithmetic
progressions of five squares over $\bQ(\sqrt{409})$, up to
equivalence, is $7^2, 13^2, 17^2, 409, 23^2$}. Therefore    $Q(5)=4$ because $Q(5)$ cannot be $5$ by Fermat's result.

We will prove that the arithmetic progression $24n+1$ is the
only one, up to equivalence, that contains $Q(N)$ squares for the
values of $N$ such that $Q(N)$ increases, for $7\le N\le 52$
 (hence, for $N=8,13,16,23,27,36,41$ and $52$). This
allow us to assert, what we have called
Super-Strong Rudin's Conjecture: let be $N=\mathcal{GP}_k+1\ge 8$
for some integer $k$, where $\mathcal{GP}_k$ is the $k$-th
generalized pentagonal number, then $Q(N)=Q(N;q,a)$ with
$\gcd(q,a)$ squarefree and $q> 0$ if and only if $(q,a)=(24,1)$.

The following theorem summarizes the main results of this article.

\begin{thm}
Let $N$ be a positive integer, then:
\begin{itemize}
\item[(S)] $Q(N)=Q(N;24,1)$ if  $6\le N\le 52$. \item[(SS)]  If
$8\le N=\mathcal{GP}_k+1\le 52$ for some integer $k$, then
$Q(N)=Q(N;q,a)$ with $\gcd(q,a)$ squarefree and $q\ge 0$ if and
only if $(q,a)=(24,1)$.
\end{itemize}
\end{thm}

\section{Preliminaries: Equivalences, translation to geometry and notations.}\label{preliminaries}
We denote by $\bN$ the set of non-negative integers.

Observe first of all, that in order to compute $Q(N)$, there is no
lost of generality to only consider  the arithmetic progressions
$qn+a$ with $\gcd(q,a)$ squarefree.

For any subset $I\subseteq \bN$, we denote by
$$\cZ_I=\{ (q,a)\in \bZ^2\ | \ \gcd(q,a)  \mbox{ squarefree },\  q\ne0 \mbox{ and } qi+a \mbox{ is a square } \forall i\in I\},$$
and by $z_I=\# \cZ_I$ its cardinality.  We will prove that if $I$ is a finite subset of $\bN$ of cardinality bigger than $3$ then $z_I<\infty$. Moreover, notice that, if $J\subseteq I$, then $\cZ_I\subseteq \cZ_J$, so $z_I\le z_J$. Since
we are interested on the subsets $I$ such that $z_I=0$, if there is
some subset $J$ with $z_J=0$, the same is true for all the subsets
$I$ containing $J$. Observe also that for all $I \subset\{0,\dots,N-1\}$ with $\# I>Q(N)$ we have by
definition that $z_I=0$.


Given $q$ and $a$ integers, $q\ne 0$, we denote by
$$\cS(q,a)=\{i\in \bN \ |\ qi+a \mbox{ is a square }\}$$
the set of values where the arithmetic progression takes its
squares, and, given $N\ge 2$, by $\cS_N(q,a)$ the set
$\cS(q,a)\bigcap \{0,1,\dots,N-1\}$. Therefore $\# \cS_N(q,a) \le Q(N)$.

\

The following lemma is elementary, and its proof is left to the reader.

\begin{lemma}
\mbox{$\quad$}
\begin{itemize}
\item[(i)] $\cS(24,1)=\left\{\mathcal{GP}_k\right\}_{k\in \bZ}$ the progression\footnote{$\mathcal{GP}_k=k(3k-1)/2$ is the sequence A001318 in \cite{oeis}.} of generalized pentagonal num\nolinebreak bers.
\item[(ii)] $\cS(8,1)=\left\{\mathcal{T}_k\right\}_{k\in \bN}$ the progression\footnote{$\mathcal{T}_k=k(k+1)/2$ is the sequence A000217 in \cite{oeis}} of triangular numbers.
\end{itemize}
\end{lemma}

Given any subset of the naturals numbers $I\subset \bN$, we will
numerate its elements by increasing order starting from $n_0$, so if
$n_i$ and $n_{i+1}$ are elements in $I$, then $n_i<n_{i+1}$ and
there is no element $m\in I$ such that $n_i<m<n_{i+1}$.

We define the following three operations on the subsets $I\subset \bN$:\\
$\bullet$ Let $i\in \bZ$ such that $n_0+i>0$, we denote by $I+i$ the translated subset
$$I+i=\{j\in \bN\ | \ j-i\in I\}.$$
$\bullet$ Let $r\in \bQ^*$ such that $ri\in \bN$ for all $i\in I$, we
denote by $rI$ the subset of $\bN$ defined by
$$rI=\{ri\ |\ i\in I\}.$$
$\bullet$ If $I$ is a finite set, $I=\{n_0,\dots,n_k\}$, we denote
by $I^s$, the symmetric of $I$, as
$$I^s=\{n_0+n_k-i\ | \ i\in I\}.$$
Therefore a finite subset is symmetric if $I^s=I$.

We say that two finite subsets $I$ and $J$ of $\bN$ are
equivalent, and denote by $I\sim J$, if there exists
$I=I_0,I_1,\dots,I_k=J$ finite subsets of $\bN$ such that
$I_{i+1}=I_i+j$ or $I_{i+1}=rI_i$ or $I_{i+1}=I_i^s$, for all
$i=1,\dots,k-1$.

Given a finite subset $I$ of $\bN$, we will denote by $n_I$ the
positive integer $\sum_{i\in I} 2^i$. Then we have a bijection
between the set of finite subsets of $\bN$ and $\bN$ (the empty set
corresponding to 0). Given  two finite subsets $I$ and $J$ of $\bN$,
we will say that $I<J$ if $n_I<n_J$.

We say that a finite subset $I$ of $\bN$ is primitive if $0\in I$,
the elements of $I$ are coprime and $n_I\le n_{I^s}$. Then any
finite subset of $\bN$ is equivalent to a unique primitive subset.

\begin{lemma} Let $I$ and $J$ be two finite subsets of $\bN$. If $I\sim J$, then $z_I=z_J$.
\end{lemma}

\begin{pf}
This is a straightforward computation, just by checking the three
possible elementary equivalences. For example, if $J=mI$, for
$m\in \bZ_{>0}$, we may assign to every element $(q,a)\in \cZ_J$,
the element $(\frac{qm}{d^2},\frac{a}{d^2})\in \cZ_I$, where $d^2$
is the largest square dividing $\gcd(qm,a)$. And to every element
$(q,a)\in \cZ_I$, the element
$(\frac{qm}{d^2},\frac{am^2}{d^2})\in \cZ_J$, where
$d=\gcd(q,a,m)$. The assertion in the case $J=I+i$ is elementary.

In the last case $J=I^s$, we may restrict ourselves to the
case $0\in I$. Then, we only need to observe that to every element
$(q,a)\in \cZ_J$, the element $(-q,a+q(N-1))$ is in $\cZ_I$.
\end{pf}

\

Let $I=\{n_0,n_1,\dots, n_k\}\subset \bN$ be a finite subset such that $k>1$ and $K$ be a field. We denote by $C_I$ the curve in $\bP^{k}(K)$ defined by the system of equations
$$
C_I\,:\,\left\{
(n_{i+2}-n_{i+1})X_{i-1}^2-(n_{i+2}-n_{i})X_{i}^2+(n_{i+1}-n_i)X_{i+1}^2=0
\right\}_{i=1,\dots, k-1}.
$$

This curve is defined over any field, and, if the characteristic of
the field is not $2$, it contains $2^{k}$ trivial points $\cT_I$
corresponding to the values $X_i^2=1$ for all $i=0,..,k$.

The following proposition collects some useful facts about the curves $C_I$ that will be used in the sequel.

\begin{prop} \label{propCI}
Let $I=\{n_0,n_1,\dots, n_k\}\subset \bN$ be a finite subset such that $k>1$ and $C_I$ be the associated curve. Then
\begin{enumerate}
\item If $K$ is a field with characteristic $0$,
then the curve $C_I$ is a non-singular projective curve of genus
$g_{k}=(k-3)2^{k-2}+1$.

\item If $k>2$, then for any $i\in I$, the natural map $C_I\to
C_{I\setminus \{i\}}$ is of degree $2$ and ramified on the $2^{k-1}$
points with $X_i=0$.

\item If $J$ is another finite set with $I\sim J$, then $C_I\cong
C_J$, with the isomorphism being the identity or the natural
involution in $\bP^{k}$ given by $$[x_0,\dots,x_{k}] \mapsto
[x_{k},\dots,x_0].$$

\item Consider the map $\iota:C_I\to \bP^{k}$ defined by
$\iota([x_0:\dots:x_k])=[x_0^2:\dots:x_k^2]$. Then there is a
natural bijection between $\iota(C_I(\bQ)\setminus \cT_I)$ and
$\cZ_I$.\
\end{enumerate}
\end{prop}

\begin{pf}
The first three items are well--known facts on the curves $C_I$ (cf. \cite{Bombieri-Granville-Pintz1992}). Now, we prove last statement. By (3), we may suppose that $I$ is
primitive, and, in particular, $n_0=0$. On one hand, let
$[x_0:\dots:x_k]\in C_I(\bQ)\setminus \cT_I$, and without lost of generality assume that $x_0,\dots,x_k$ are coprime integers. Then the corresponding element of $\cZ_I$ is
$((x_1^2-x_0^2)/n_2,x_0^2)$. On the other hand, let
$(q,a)\in \cZ_I$. Then the point $[a:a+n_1q:\dots:a+n_kq]$ belongs to $\iota(C_I(\bQ)\setminus \cT_I)$.
\end{pf}

\

In fact, in order to compute $\cZ_I$, we may carry out  it
modulo the group of automorphisms $\Upsilon_I$ of the curve $C_I$
generated by the automorphisms $\tau_i(x_i)=-x_i$ and
$\tau_i(x_j)=x_j$ if $i\ne j$, for $i=0,\dots,k$. Notice that
$C_I(\bQ)/\Upsilon\cong \im(\iota)$.

\

Finally, observe that, if $I\subset\bN$ only has three elements, the corresponding curve $C_I$ is a genus 0 curve. And since it has rational points, $C_I(\bQ)=\bP^1(\bQ)$ and hence
$z_I=\infty$.

\section{First elementary cases}\label{first_cases}

In order to study the first values of $N$, we will consider the subsets $I\subset \bN$ of cardinality $4$. In this case the curve $C_I$ has genus $1$. Moreover,  $C_I$ is an elliptic curve over any field, since $C_I$ always has the rational points $\cT_I$.

\begin{prop}\label{firstcases} Given any subset $I$ of $\{0,1,\dots,6\}$ with four
elements, we have that $z_I=\infty$ unless $I$ is equivalent to one
of the following five subsets:
$$ \{0,1,2,3\},\ \{0,1,3,4\},\ \{0,1,4,5\}, \ \{0,2,3,5\} \mbox{ and } \{0,1,5,6\}.$$
In which case, $z_I=0$.
\end{prop}

\begin{pf} To prove that $z_I=0$ for the given finite sets $I$ in the
proposition, one only needs to show that $C_I(\bQ)=\cT_I$. Or, equivalently, that $\# C_I(\bQ)=8$. Using standard
transformations (see section \ref{4squares}), one may put the corresponding
elliptic curves in Weierstrass form, and then compute the Cremona
reference. We obtain the elliptic curves \verb+24a1+,
\verb+48a3+, \verb+15a3+, \verb+120a2+ and \verb+240a3+ respectively. All of them have rank $0$ and
torsion subgroup isomorphic to $\bZ/2\bZ\oplus \bZ/4\bZ$. Hence $\#
C_I(\bQ)=8$ in these cases. For all the other cases, one easily shows 
that the rank of the corresponding elliptic curve $C_I$ is 1. Therefore, they have an infinite number of points.
\end{pf}

\

Observe that all the subsets in the proposition with $z_I=0$ are
symmetric subsets. We will see that this is true for any subset with
four elements.

\begin{cor}
$Q(6)=Q(7)=4$ and $Q(8)=5$.
\end{cor}

\begin{pf} First of all, we clearly have that $4=Q(5)\le Q(6)\le 5$.
Suppose we have a primitive subset $I$ of $5$ elements inside
$\{0,1,2,3,4,5\}$ such that $z_I>0$. Notice that if $I$ does
not contain $5$, $I\setminus\{5\}$ will be $\{0,1,2,3,4\}$, which
contains $J=\{0,1,2,3\}$. Since $z_J=0$ by the proposition
\ref{firstcases}, $z_I=0$. So $I=\{0,1,2,3,4,5\}\setminus\{i\}$ for $0<i<5$, and
we have 4 cases: if $i=1$ or $i=4$, then $I$ contains $\{2,3,4,5\}\sim \{0,1,2,3\}$
or $\{0,1,2,3\}$ respectively, and we have again the same result; and if $i=2$ or $i=3$, then $I$ contains $\{0,1,3,4\}$ or $\{1,2,4,5\}\sim \{0,1,3,4\}$, and applying another case of proposition \ref{firstcases}, we conclude.

Now, we are going to prove that $Q(7)=4$. We again proceed with the same strategy. We
consider $I$ a primitive subset of $\{0,\dots,6\}$ with five
elements and we show that $z_I=0$ if we find a subset $J$ of $I$
with four elements appearing in the list of the proposition
\ref{firstcases}, hence with $z_J=0$. We may suppose that
$I=\{0,i,j,k,6\}$ for some $0<i<j<k<6$. We have $10$ cases, but only $6$ cases to consider because the symmetries: the first case $\{0,1,2,3,6\}$ is by Fermat; the second case $\{0,1,2,4,6\}$ because it contains $\{0,2,4,6\}\sim \{0,1,2,3\}$, hence again by Fermat; $\{0,1,2,5,6\}$ contains $\{0,1,5,6\}$, $\{0,1,3,4,6\}$ contains $\{0,1,3,4\}$, $\{0,1,3,5,6\}$ contains $\{0,1,5,6\}$ and the last subset $\{0,2,3,4,6\}$ contains $\{0,2,4,6\}\sim \{0,1,2,3\}$.

Now, since $Q(8)\le Q(7)+1=5$, to prove $Q(8)=5$ we only need to exhibit
an arithmetic progression with five squares in the first 8 terms,
and the arithmetic progression $1+24n$ do the job. Note that
$\cS_8(24,1)=\{0,1,2,5,7\}$.
\end{pf}

\

So, the first strategy to detect subsets $I$ with cardinality bigger that $3$ and
$z_I=0$ is to find a subset $J$ of $I$ with four elements such that
$z_J=0$. This will be carry out considering the associated elliptic curve $E_J$ and
showing that it only contains the (eight) trivial points $\cT_J$. In
the next section we will study some necessary conditions where this happens.

\section{Four squares in arithmetic progressions}\label{4squares}

Let be $I=\{n_0,n_1,n_2,n_3\}$ with $n_0<n_1<n_2<n_3$. For
any $(q,a)\in \cZ_I$, consider the eight points in the three
dimensional projective space $[x_0,x_1,x_2,x_3]\in \bP^3$ such
that $x_j^2=qn_{j}+a$. They all lie in the curve $C_I$ given by
the sytem of equations
$$
C_I\,:\,\left\{
\begin{array}{rcl}
(n_2-n_1)X_0^2-(n_2-n_0)X_1^2+(n_1-n_0)X_2^2=0,\\[1mm]
(n_3-n_2)X_1^2-(n_3-n_1)X_2^2+(n_2-n_1)X_3^2=0.
\end{array}
\right.
$$
We have that $C_I$ is an elliptic curve since it has genus 1 (being the intersection of two
quadric surfaces in $\bP^3$) and it has the rational points $[1,\pm1,\pm1,\pm1]$.

Let us denote
$$m_0=\frac{n_1-n_0}{n_2-n_1}\,,\quad \ m_1=\frac{n_3-n_2}{n_2-n_1}.$$
Note that they are both strictly positive rational numbers. Then we may
write the equations of $C_I$ as
$$
C_I\,:\,\left\{
\begin{array}{rcl}
X_0^2-(m_0+1)X_1^2+m_0X_2^2=0,\\[1mm]
m_1X_1^2-(m_1+1)X_2^2+X_3^2=0.
\end{array}
\right.
$$
Now, we parametrize the first equation as 
$$
[X_0:X_1:X_2]=
[(m_0+1)-2(m_0+1)t+t^2\,:\,(m_0+1)-2t+t^2 \,:\,(m_0+1)-t^2].
$$
Next, we substitute $X_0,X_1,X_2$ in the second equation and we obtain a new equation
of the curve, depending on the parameter $t$ (note that $t=(X_2-X_0)/(X_2+ X_1)$):
$$
C_I: X_3^2=t^4+4m_1t^3-2(m_0+4m_1+2m_1 m_0+1)t^2+4m_1(m_0+1)t+(m_0+1)^2.
$$
Observe that the set $I$ is symmetric if and only if $m_0=m_1$.

\begin{lemma} \label{EI} Let $E_I$ be the elliptic curve defined by the Weierstrass form
$$ E_I\,:\,y^2=x(x-m_0m_1)(x+m_0+m_1+1).$$
Then, there exists a $\bQ$-isomorphism $\phi:C_I\longrightarrow E_I$ such that $\phi([1,1,1,1])=[0,1,0]$. Furthermore, if we denote by $\cF_I=\phi(\cT_I)$, then $\#\cF_I=8$.
\end{lemma}

\begin{pf}
The proof of the existence of the isomorphism $\phi$ is an straightforward computation. For example, it may be carried out using the formulae in \cite[section 7.2]{Cohen239}. The set $\cF_I$ is described by the table below
$$
\begin{array}{|c|c|cc|}
\hline
i & P_i & & Q_i=\phi(P_i) \\[1mm]
\hline
0 & [1,1,1,1] & & \mathcal O=[0,1,0] \\[1mm]
\hline
1 & [-1,1,-1,1] & & (0,0)\\[1mm]
\hline
2 & [-1,1,1,-1] & & (m_0 m_1,0)\\[1mm]
\hline
3 & [-1,-1,1,1] & & (-m_0-m_1-1,0)\\[1mm]
\hline
4 & [1,1,-1,1] & & (-m_1,-m_1(m_0+1))\\[1mm]
\hline
5 & [1,-1,1,1] & & (-m_0,m_0(m_1+1))\\[1mm]
\hline
6 & [-1,1,1,1] & & (m_0(m_0+m_1+1),-m_0(m_0+1)(m_0+m_1+1)) \\[1mm]
\hline
7 & [1,1,1,-1] & & (m_1(m_0+m_1+1),m_1(m_1+1)(m_0+m_1+1))\\[1mm]
\hline
\end{array}
$$
Therefore, we have $\#\cF_I=8$ since $m_0,m_1>0$.
\end{pf}

\begin{cor} The set $\cF_I$ is a subgroup of $E_I$ if and only if
$I$ is symmetric. Furthermore, if $I$ is not symmetric, then $z_I>0$.
\end{cor}

\begin{pf} First, observe that the opposites of the non-Weierstrass points on $\cF_I$
do not belong to $\cF_I$ unless $m_0=m_1$.  Since $m_0$ and $m_1$ are strictly positive numbers, then, for example, $-Q_5=\left(-m_0,-m_0(m_1+1)\right)\in\cF_I$ if and only if it is equal to
$Q_4$. This shows one implication.
Finally, if $m_0=m_1$, one easily checks that the non-Weierstrass
points are of order 4, and their doubles are equal to the point
$(m_0m_1,0)$. That is, $\cF_I\cong \bZ/2\bZ\oplus \bZ/4\bZ$.
\end{pf}

\

\begin{rem}\label{ap567} One may use the isomorphism $\phi$ in order to find explicitly which arithmetic progression
corresponds to the set of points $\{-P_4,-P_5,-P_6,-P_7\}$. If we
suppose $I$ is primitive, in particular with $n_0=0$, so it is of
the form $\{0,n_1,n_2,n_3\}$, with $n_1$, $n_2$ and $n_3$ coprime,
then the arithmetic progression $a+nq$ given by
$$
\left\{
\begin{array}{l}
a =((n_1+n_2-n_3)^2-4n_1n_2)^2,\\
q=2^3(n_1+n_2-n_3)(n_1-n_2-n_3)(n_1-n_2+n_3),
\end{array}
\right.
$$
has squares for $n=0,n_1,n_2,n_3$. Using this construction, we show in the table below the arithmetic progression $(q,a)$ for all the equivalence classes of $4$-tuples $I\subset \{0,\dots,N-1\}$, for $5\le N\le 7$, such that $C_I(\bQ)\ne \cT_I$. Note that in all of these cases $\mbox{rank}\,E_I(\bQ)=1$.\\
 \begin{center}
\begin{tabular}{|c|c|c|}
\hline
$N$ & $I$ & An Arithmetic progression $(q,a)$ such that $I\subset\mathcal{S}(q,a)$\\
\hline
$5$ & $\{0,1,2,4\} $ & $(120,49)$\\
\hline
\multirow{2}{*}{$6$}  & $\{0,1,2,5\}$ & $(24,1)$\\
                                 &  $\{0,1,3,5\}$ & $(168,121)$\\
\hline
\multirow{4}{*}{$7$}  & $\{0,1,2,6\}$ & $(840,1)$\\
                                 &  $\{0,1,3,6\}$ & $(8,1)$\\
                                 &  $\{0,2,3,6\}$ & $(280,529)$\\
                                 &  $\{0,1,4,6\}$ & $(24,25)$\\
                                 \hline
\end{tabular}
\end{center}

\end{rem}

\begin{rem} One even show that, if $I$ is not symmetric, the
subgroup generated by $\cF_I$ is infinite, unless $m_1^2+m_1+1$ is
a square and $m_0=-(m_1+2-2\sqrt{m_1^2+m_1+1})/3$, or $m_1^2+m_1$ is a square and
$m_0=m_1+2\sqrt{m_1^2+m_1}$, which are in fact the same case by
interchanging $m_0$ and $m_1$. Or equivalentment, $E_I$ is
$\bQ$-isomorphic to
$$
E_{t}\,:\,y^2= x(x+1+4t)(x+16t^3(t+1)),\qquad \mbox{for some $t\in\bQ$}.
$$
In this case,
$$
\langle P \,:\, P\in \cF_I\rangle =\cF_I\cup\{-P_4,-P_5,-P_6,-P_7\}\cong \bZ/2\bZ\oplus \bZ/6\bZ.
$$
Furthermore, the points $-P_4,-P_5,-P_6,-P_7$ correspond to the arithmetic progression with $a=0$ and $q=1$. Thus, if we suppose that $I$ is primitive then there exist $s_1,s_2\in \bZ$ such that $n_1=s_1^2$, $n_2=s_2^2$ and $n_3=(s_1\pm s_2)^2$.
\end{rem}

We have seen that if $I\subset \bN$ has four elements, a necessary condition to have $z_I=0$ is that $I$ is
symmetric. In the sequel, we obtain more necessary conditions, some of them under the Parity Conjecture.

\

Observe that the number of symmetric subsets with
$4$ elements contained in $\{0,1,\dots, N\}$ may be explicitly
computed in terms of $N$ (and it is the sequence A002623 in
\cite{oeis}, with $n=N-3$), and it is almost equal to a polynomial
of degree $3$ in $N$:
$$\frac
{N^3}{12} -\frac{7N^2}{8}+\frac{35N}{12} -
\frac{49}{16}+\frac{(-1)^N}{16}.$$ Since the number of subsets
with 4 elements is a polynomial of degree $4$ in $N$, there are
$2/N+O(N^{-2})$ symmetric subsets among all the subsets with four
elements. We do not know how many of the equivalence clases of
subsets with four elements are symmetric, but we suspect is of the
same order.

\

In the primitive symmetric case, that is if
$I=\{0,n_1,n_2,n_1+n_2\}$ with $0<n_1<n_2$ coprime integers, we
have that a $\bQ$-isomorphism $\psi:C_I\longrightarrow E'_t$ exists, where
$$
E'_{t}\,:\,y^2=x(x+1)(x+t^2),\qquad \mbox{and $t=n_2/n_1$}.
$$
In this case, $\psi(\cT_I)=\{\cO,(0,0),(-1,0),(-t^2,0),(\pm t, \pm t(t+1))\}$.

\

The first remark is that there are plenty of symmetric
subsets $I$ with four elements such that $z_I>0$, and even with
infinite number of elements. For example, this is the case if the torsion subgroup has more than $8$ elements. This case only may occur, by Mazur's theorem, if the torsion subgroup of $E_{t}(\bQ)$ is isomorphic to
$\bZ/2\bZ\oplus \bZ/8\bZ$. Or, equivalently, some of the four
$4$-torsion points given by the points with  $x$-coordinate equal to $\pm t$ is the
double of some rational point. We use the standard formulae (or
by a 2-descent argument) to obtain that this happens if and only if
the $x$-coordinate and the $x$-coordinate $+1$ are squares.

\begin{lemma} Let $t$ be a positive rational number, then
$$
E'_t(\bQ)_{tors}\cong
\left\{
\begin{array}{ccl}
\bZ/2\bZ\oplus \bZ/8\bZ & & \mbox{if $t=\left(s-\frac 1{4s}\right)^2$, for some $s\in\bQ$,}\\[1mm]
\bZ/2\bZ\oplus \bZ/4\bZ & &\mbox{otherwise.}
\end{array}
\right.
$$
\end{lemma}

In fact, we may exactly characterize which primitive sets have
their corresponding elliptic curve with torsion subgroup of order
$16$, and even which arithmetic progressions correspond to the
torsion points.

\begin{cor} Let $I=\{0,n_1,n_2,n_1+n_2\}$ be a primitive symmetric subset of $\bN$.
Then the elliptic curve $E_I$ has torsion isomorphic to
$\bZ/2\bZ\oplus \bZ/8\bZ$ if and only if $n_1$, $n_2$ and $n_1+n_2$
are squares. The torsion points of $E_I$ correspond to the constant
arithmetic progression along with the arithmetic progressions
with $a=0$ and $q=1$ and with $a=n_1+n_2$ and $q=-1$.
\end{cor}

\begin{pf} Since $t=n_2/n_1$ is a square, both $n_1$ and $n_2$
are squares. Since $t+1=(n_1+n_2)/n_1$ is also a square, also
$n_3=n_1+n_2$ is a square. Trivially, the arithmetic progression
with $a=0$ and $q=1$ verifies that $a_{n_i}=a+n_iq=n_i$ are squares,
and the one with $a=n_1+n_2$ and $q=-1$ verifies that
$a_{n_i}=a+n_iq=n_{3-i}$ are also squares. And they correspond to the torsion points of order $8$.
\end{pf}

\

Concerning the general symmetric case, we find one parametric subfamilies, and even two parametric subfamilies of $E'_t$, with other rational points apart of the trivial ones.

\begin{exmp} Let $z_1$ and $z_2$ be non-zero rational numbers, and
consider
$$t=\frac 14 \left( z_1+\frac 1{z_1} + z_2+\frac 1{z_2}\right)\mbox{ and }
x=\frac{-(z_1+z_2)^2}{4z_1z_2}.
$$
Then $x(x+1)(x+t^2)$ is a square. Moreover, if $z_1,z_2\ne
\pm1$,  $z_1 \ne \pm z_2$ and $z_1 \ne 1/z_2$, then we obtain a non-trivial point in $E'_t$.
\end{exmp}

Recall that the Parity Conjecture asserts that any elliptic curve $E$
defined over $\bQ$ (or a general number field) that has root
number $W(E)=-1$ has odd rank (and, in particular, infinitely many
rational points). But the root number is easily computable, even
in some families. In our case we have an explicit description.

\begin{prop} Let $0<a<b$ be coprime integers, and let $E_{a,b}$ be the elliptic curve defined by
$$
E_{a,b}\,:\,y^2=x(x+a^2)(x+b^2).
$$
Then $W(E)=-1$ if and only if $\alpha(a,b)  \equiv \mu_2(a,b) \pmodd{2}$, where
$$
\begin{array}{l}
\alpha(a,b)=\#\{p \mbox{ odd prime} \, | \, \mbox{$p$ divides
$ab$}\} \\
\qquad\qquad\qquad\qquad\qquad\quad +\, \#\{p \mbox{ prime}
\, | \, \mbox{$p$ divides $(b^2-a^2)$ and  $p\equiv 1 \pmodd{4}$} \},\\[2mm]
\mu_2(a,b)=\left\{
\begin{array}{ll} 0 & \mbox{ if } ab \equiv 4 \pmodd{8} \mbox{ or } ab\equiv
1 \pmodd{2} \mbox{ and } \frac{(b^2-a^2)}{8} \equiv 1 \pmodd{2},
\\1 & \mbox{ otherwise.}
\end{array}\right.
\end{array}
$$
\end{prop}

\begin{pf} Recall that the root number $W(E)$ of an elliptic curve $E$ over
$\bQ$ is equal to the product of local root numbers $W_p(E)$, where $p$ runs over all the prime numbers and infinity. One always has that $W_{\infty}(E)=-1$, $W_p(E)=1$ if the curve has
good or non-split multiplicative reduction  at $p$, and
$W_p(E)=-1$ if the curve has split multiplicative reduction at
$p$. Since $a$ and $b$ are coprime integers, $E_{a,b}$ is minimal at an odd prime. The reduction is
good if $p$ does not divide $ab(b^2-a^2)$, and it is split
multiplicative if $p$ divides $ab$, or if $p$ divides $b^2-a^2$
and $-1$ is a square modulo $p$. Hence we obtain $W(E)=(-1)^\beta W_2(E)$, where $\beta=\alpha(a,b)+1$.

Now, to compute the root number at 2, we need to carry out a more
detailed analysis, since the reduction can be additive in this
case. First, one shows by a change of variables that the curve
$E_{a,b}$ is isomorphic to the curve given by the equation
$y^2+xy+ry=x^3+rx^2$ (this is the so called Tate normal form),
where $r=\frac{ab}{4(a+b)^2}$. Given $s\in \bQ$, we denote by
$v_2(s)$ the $2$-adic valuation of $s$. The curve has good
reduction at $2$ (hence $W_2(E)=1$) if and only if $v_2(r)=0$, so
if $v_2(a)=2$ or $v_2(b)=2$. And it has split multiplicative reduction
(hence $W_2(E)=-1$) if $v_2(r)>0$, so if $v_2(a)>2$ or $v_2(b)>2$,
and it has additive reduction otherwise.

When the valuation of $r$ is negative, we will consider the
original equation of the curve $E_{a,b}$, which is an integral
model. We will need the following standard invariants of the
equation:
$$
\begin{array}{lccl}
\displaystyle j=j(E_{a,b})=\frac{2^8(a^4+b^4-a^2b^2)^3}{a^4b^4(a+b)^2(a-b)^2}, & & & c_4=2^4(a^4+b^4-a^2b^2),\\[4mm]
\Delta=2^4a^4b^4(a+b)^2(a-b)^2,& & & c_6=2^5(a^2+b^2)(2b^2-a^2)(2a^2-b^2).
\end{array}
$$

Now, if $v_2(r)=-1$, hence $v_2(a)=1$ or $v(b)=1$. Then the curve
has potentially good reduction ( since $v_2(j)=4$), and we can look
at the tables in \cite{Halberstadt1998}. We obtain $v_2(\Delta)=8$,
$v_2(c_4)=4$ and $v_2(c_6)=6$, hence this information is not
enough to obtain the sign. We consider $c_6'=c_6/2^6$ and
$c_4'=c_4/2^4$, and we compute $2c_6'+c_4' \pmodd{16}$. After a
case by case computation one obtain that $2c_6'+c_4'\equiv 7
\pmodd{16}$, so we are in the case $I_1^*$. Next, one
computes that $2c_6'+c_4' \equiv 23 \not{\!\!\equiv} 7
\pmodd{32}$, so the root number $W_2(E)$ is $-1$ in this case.

Now, if $v_2(a)=v_2(b)=0$, we need to take in account the
valuations $v_2(a+b)$ and $v_2(a-b)$. First we need to determine
when the curve $E_{a,b}$ has potentially multiplicative reduction
in order to apply the formulae by Rohrlich \cite{Rohrlich1993}. This
is equivalent to the case $v_2(j)<0$. Since in our case
$v_2(j)=8-2v_2(a-b)-2v_2(a+b)$, we have potentially
multiplicative reduction if and only if $v_2(a-b)+v_2(a+b)\ge 4$.
In this case, the root number is computed as follows: if $s\in \bQ$, we denote by $\overline{s}=s2^{-v_2(s)}$. Then, we obtain that $W_2(E)\equiv -\overline{c_6} \pmodd{4}$. But observe that
$$\overline{c_6}=\overline{(a^2+b^2)}(2b^2-a^2)(2a^2-b^2)\equiv
\overline{(a^2+b^2)} \equiv 1 \pmodd{4},$$ hence $W_2(E)=-1$ in
this case.

Finally, we need to consider the potentially good reduction case,
with $v_2(a)=v_2(b)=0$. In this case, we have
$v_2(a-b)+v_2(a+b)=3$. Therefore, we obtain that $v_2(\Delta)=10$,
$v_2(c_4)=4$ and $v_2(c_6)=6$. One easily shows that the necessary
condition $\overline{c_6} \equiv 1 \pmodd{4}$ for the case $I_2^*$
in the tables of \cite{Halberstadt1998} is always satisfied, so we obtain that $W_2(E)$ is $1$ in this case.

In summary, we obtain $W_2(E)=1$ if and only if $v_2(ab)=0$ and
$v_2(a^2-b^2)=3$, or $v_2(ab)=2$.
\end{pf}

\

\begin{cor} Let $0<n_1<n_2$ be coprime integers and $I=\{0,n_1,n_2,n_1+n_2\}$. Assume
that the Parity Conjecture holds for the curve $E_{n_1,n_2}$, and
that $\alpha(n_1,n_2)  \equiv \mu_2(n_1,n_2) \pmodd{2}$. Then $z_I=\infty$.
\end{cor}

\begin{rem}
Cohn \cite{Cohn1983} studied the special symmetric case $\{0,2,n,n+2\}$ when $n\leq 100$. The corollary above gives an arithmetic sufficient condition to determine if there is an arithmetic progression with squares at the positions $\{0,2,n,n+2\}$ for any positive integer $n$. This condition is that the number of odd primes dividing $n$ has the same parity than the number of primes congruent to $1$ mod $4$ dividing $n^2-4$. The disadvantage is that this condition is under the Parity Conjecture for the elliptic curve $E_{2,n}$. Note that we may suppose that $n$ is odd since in the even case we can reduce $\{0,2,n,n+2\}$ to $\{0,1,n/2,(n+1)/2\}$.
\end{rem}

\section{Five squares in arithmetic progressions: the technique}\label{5squares}

We will see in section \ref{sec_comp_result} that the results of section \ref{4squares} are not
enough to show the Rudin's Conjecture even for small values of N. In
this section we will study how to prove, for some subsets $I$ with
5 elements, that $z_I=0$ even if it is not zero for any subset $J$ of
$I$ with 4 elements. Moreover, we will be able to determine $\cZ_I$ in some cases.

In order to prove these type of results, we need to be able to
compute the rational points of some genus 5 curves whose Jacobians
are product of elliptic curves, all of them of rank greater than
0. Hence it is not possible to apply the classical Chabauty
method (see \cite{Chabauty1941}, \cite{Coleman1985}, \cite{Flynn1997}, \cite{Stoll2006}, \cite{Stoll2007}, \cite{McCallum-Poonen}). We will
instead apply the covering collections technique, as developed by
Coombes and Grant \cite{Coombes-Grant1989}, Wetherell \cite{Wetherell1997} and others, and
specifically a modification of what is now called the elliptic curve
Chabauty method developed by  Flynn and Wetherell in \cite{Flynn-Wetherell2001} and
by Bruin in \cite{Bruin2003}. In
fact, we will follow the same technique we applied in \cite{Gonzalez-Jimenez-Xarles2011},
though we could may be use also a similar technique as the one we
use in \cite{Gonzalez-Jimenez-Xarles2013} to study 5 squares in arithmetic progression
over quadratic fields.

First, we fix the notations. We consider a primitive subset $I\subset \bN$ with 5 elements, and $C_I$ the associated curve as in
proposition \ref{propCI}. Then if we want to prove that $z_I=0$, this will be equivalent to prove that
$C_I(\bQ)$ only contains the trivial points $\cT_I$. In fact,
since the genus of $C_I$ is greater than one, its set of rational
points always is finite, hence we may even try to compute them.

Observe that $C_I$ has 5 different maps to the elliptic curves
corresponding to $C_J$, for $J$ a subset of $I$ with four
elements. As we have already seen in the previous section, the
corresponding elliptic curves have all its $2$-torsion points
defined over $\bQ$, a fact that we will use to build
unramified coverings of $C_I$.

The method has two parts. Suppose we
have a curve $C$ over a number field $K$, and an unramified map
$\chi:C'\to C$ of degree greater than one, may be defined over a
finite extension $L$ of $K$, along with a nice quotient $C'\to
H$, for example a genus $1$ quotient. We consider the different
unramified coverings $\chi^{(s)}:C'^{(s)}\to C$ build by all the
twists of the given one. We obtain that
$$
C(K)=\bigcup_{s} \chi^{(s)}(\{P\in C'^{(s)}(L)\ : \ \chi^{(s)}(P)\in C(K)\}) ,
$$
the union being disjoint. Moreover, only a finite number of twists
have rational points, and the finite set of twists
with points locally everywhere can be explicitly described. The
method depends first on being able to compute the set of twists,
and second, on being able to compute the points $ P\in
C'^{(s)}(L)$ such that $\chi^{(s)}(P)\in C(K)$, by computing their images in $H^{(s)}(L)$.
$$
\begin{array}{ccc}
\xymatrix@R=0.7pc@C=0.8pc{
 C' \ar@{->}[dd]_\chi  \ar@/^2mm/[ddrr]^\pi& & \\
     & &\\
 C & & H
}
&
\xymatrix@R=0.7pc@C=0.8pc{
 & & & &\\
  \ar@{~>}[rrrr]^s & & & &\\
 & & & &
}
&
\xymatrix@R=0.7pc@C=0.8pc{
 C^{'(s)} \ar@{->}[dd]_{\chi^{(s)}}  \ar@/^2mm/[ddrr]^{\pi^{(s)}}& & \\
     & &\\
 C & & H^{(s)}
}
\end{array}
$$

In our case, the coverings we are
searching for will be defined over $\bQ$, but the genus $1$
quotients of such coverings are, in general, not
defined over $\bQ$, but in a quadratic or in a biquadratic
extension. The way we will construct the coverings (factorizing quartic polynomials) will also give us the genus $1$ quotients and the field where they are all defined.

In order now to construct the coverings of the curve $C_I$, we
first rewrite the curve as the projectivization (and
desigularization) of a curve in $\bA^3$ given by equations of the
form $y_1^2=p_1(x)$ and $y_2^2=p_2(x)$, where $p_1(x)$ and
$p_2(x)$ are separable degree 4 polynomials with coefficients in
$\bQ$. This is possible because of the special form of the curve
(essentially, because it has two degree 2 maps to elliptic curves
that correspond to involutions that commute with each other), and
in our case we will see it can be done in 10 different ways.

Next, we will consider a factorization of the polynomials
$p_i(x)$  as product of two degree two polynomials $p_{i,1}(x)$
and $p_{i,2}(x)$, may be defined over a larger field $K$ (in our
case, a quadratic field). This factorization
$p_i(x)=p_{i,1}(x)p_{i,2}(x)$ determines an unramified degree two
covering $\chi:F'_i\to F_i$  of the genus $1$ curve $F_i$ given by
$y_i^2=p_i(x)$, as we describe in the next proposition.

\begin{prop}\label{2cov}
Let $F$ be a genus $1$ curve over a number field $K$ given by a
quartic model of the form $y^2=q(x)$, where $q(x)$ is a degree
four monic polynomial in $K[x]$. Thus, the curve $F$ has two
rational points at infinity, and we fix an isomorphism from $F$ to
its Jacobian $E=\Jac(F)$ defined by sending one of these points
at infinity to the zero point of $E$. Then
\begin{enumerate}
\item Any $2$-torsion point of the curve $E$ defined over $K$
corresponds to a factorization of the polynomial $q(x)$ as a
product of two quadratic polynomials $q_1(x),q_2(x) \in L[x]$, where $L/K$ is an algebraic extension of degree at most $2$.

\item Given such a $2$-torsion point $P$, the degree two
unramified covering $\chi:F'\rightarrow F$ corresponding to the degree two
isogeny $\phi:E'\to E$ determined by $P$ can be described as the
map from the curve $F'$ defined over $L$, with affine part in
$\bA^3$ given by the equations $y_1^2=q_1(x)$ and
$y_2^2=q_2(x)$ and the map given by $\chi((x,y_1,y_2))=(x,y_1y_2)$.

\item Given any degree two isogeny $\phi:E'\to E$, consider the
Selmer group $\Sel(\phi)$ as a subgroup of $K^*/(K^*)^2$. Let
$\cS_L(\phi)$ be a set of representatives in $L$ of the image
of $\Sel(\phi)$ in $L^*/(L^*)^2$ via the natural map. For any
$\delta\in \cS_L(\phi)$, define the curve $F'^{(\delta)}$ given
by the equations $\delta y_1^2=q_1(x)$ and $\delta y_2^2=q_2(x)$,
and the map to $F$ defined by
$\chi^{(\delta)}(x,y_1,y_2)=(x,y_1y_2/\delta^2)$. Then
$$F(K) \subseteq \bigcup_{\delta \in \cS_L(\phi)} \chi^{(\delta)}(\{(x,y_1,y_2)\in F'^{(\delta)}(L)\ : \ x\in K  \mbox{ or } x=\infty \}).$$
\end{enumerate}
\end{prop}

\begin{pf}
First we prove (1) and (2). Suppose we have such a factorization
$q(x)=q_1(x)\,q_2(x)$ over some extension $L/K$, with $q_1(x)$ and
$q_2(x)$ monic quadratic polynomials. Then the covering
$\chi:F'\rightarrow F$ from the curve $F'$ defined over $L$, with
affine part in $\bA^3$ given by the equations $y_1^2=q_1(x)$ and $y_2^2=q_2(x)$ and the map given by
$\chi((x,y_1,y_2))=(x,y_1y_2)$, is an unramified degree two
covering. So $F'$ is a genus 1 curve, and clearly it contains the
preimage of the two points at infinity, which are rational over
$L$, hence it is isomorphic to an elliptic curve $E'$. Choosing
such isomorphism by sending one of the preimages of the fixed
point at infinity to $\cO$, we obtain a degree two isogeny $E'\to E$,
which corresponds to a choice of a two torsion point.

So, if the polynomial $q(x)$ decomposes completely in $K$, the
assertions (1) and (2) are clear since the number
of decompositions $q(x)=q_1(x)\,q_2(x)$ as above is equal to the
number of points of exact order 2. Now the general case is proved
by Galois descent: a two torsion point $P$ of $E$ is defined over
$K$ if and only if the degree two isogeny $E'\to E$ is defined
over $K$, so if and only if the corresponding curve $F'$ is
defined over $K$. Hence the polynomials $q_1(x)$ and $q_2(x)$
should be defined over an extension of $L$ of degree $\le 2$, and in
case they are not defined over $K$, the polynomials $q_1(x)$ and
$q_2(x)$ should be Galois conjugate over $K$.

Now we show the last assertion. First, notice that the curves
$F'^{(\delta)}$ are twisted forms (or principal homogeneous
spaces) of $F'$, and it becomes isomorphic to $F'$ over the
quadratic extension of $L$ adjoining the square root of $\delta$.

Consider the case where $L=K$. So $F'$ is defined over $K$. For any
$\delta\in \Sel(\phi)$, consider the associated homogeneous space
$D^{(\delta)}$; it is a curve of genus $1$ along with a degree
$2$ map $\phi^{(\delta)}$ to $E$, without points in any local
completion, and isomorphic to $E'$ (and compatible with $\phi$)
over the quadratic extension $K(\sqrt{d})$. Moreover, it is
determined by such properties (see section 8.2 in \cite{Cohen239}). So,
by this uniqueness, it must be isomorphic to $F'^{(\delta)}$
along with $\chi^{(\delta)}$. The last assertion also is clear
from the definition of the Selmer group.

Now, the case $L\ne K$. The assertion is proved just
observing that the commutativity of the diagram
$$\xymatrix{
  \Sel(\phi) \ar@{->}[r]    \ar@{->}[d]& \Sel(\phi_L)   \ar@{->}[d]\\
   K^*/(K^*)^2 \ar@{->}[r]& L^*/(L^*)^2 
  }
$$
where the map $\Sel(\phi)\to \Sel(\phi_L)$ is the one sending the
corresponding homogeneous space to its base change to $L$.
\end{pf}

\

In order to apply the method to a $5$-tuple $I\subset \bN$, we first explain how to construct
models of $C_I$ as the ones described above. We need first to
choose a subset $J=\{n_0,n_1,n_2\}\subset I$ with three elements,
which determines a partition $I=J \sqcup\{n_3,n_4\}$ of $I$, the
$n_i$ not necessarily ordered, and everything will depend of that
choice. Second, we write the equations of $C_I$ of the form
$$
C_I\,:\,\left\{
\begin{array}{rcl}
X_0^2&\!\!\!\!=&\!\!\!\!(m_0+1)X_1^2-m_0X_2^2,\\[1mm]
X_3^2&\!\!\!\!=&\!\!\!\!-m_1X_1^2+(m_1+1)X_2^2,\\[1mm]
X_4^2&\!\!\!\!=&\!\!\!\!-m_2X_1^2+(m_2+1)X_2^2,
\end{array}
\right.
$$
where
$$m_0=\frac{n_1-n_0}{n_2-n_1},\qquad m_1=\frac{n_{3}-n_2}{n_2-n_1},\qquad m_2=\frac{n_{4}-n_2}{n_2-n_1}.$$
Now, we
parametrize the first equation as it has been done on section \ref{4squares}:
$$
[X_0:X_1:X_2]= [(m_0+1)-2(m_0+1)t+t^2\,:\,(m_0+1)-2t+t^2 \,:\,(m_0+1)-t^2],
$$
and 
substituting in the next two equations, we obtain the new equations
of the curve, depending on the parameter $t$:
$$
C_I:\left\{ y_1^2=p_1(t)\,,\,y_2^2=p_2(t)\right\},
$$
where, for $i=1,2$,  $y_i=X_{2+i}$ and 
$$
p_i(t)=t^4+4m_it^3-2(m_0+4m_i+2m_im_0+1)t^2+4m_i(m_0+1)t+(m_0+1)^2\,.
$$
For $i=1,2$, by lemma \ref{EI}, we obtain that the genus $1$ curve $F_i: y_i^2=p_i(t)$ is $\bQ$-isomorphic to the elliptic curve $$E_i:
y^2=x(x-m_0m_i)(x+m_0+m_i+1).$$

Next, we need to choose factorizations of the polynomials $p_i(t)$
as product of two quadratic polynomials over some quadratic
extension $K/\bQ$. We describe in the next elementary lemma all
these factorizations, relating them to the corresponding $2$-torsion points in the corresponding elliptic curve $E_i$.

\begin{lemma}\label{factorizations} For $i=1,2$, denote by
$$D_{i,1}=m_i(1 + m_i),\qquad D_{i,2}=(1 + m_i)(m_i + m_0 + 1), \qquad D_{i,3}=m_i(m_i + m_0 + 1),$$ and choose an square root $\alpha_{i,j}=\sqrt{D_{i,j}}$.
Then the polynomial $p_i(t)$ factorizes over $\bQ(\alpha_{i,j})$ as a product of two quadratic
polynomials $p_{i,j,+}(t)$ and $p_{i,j,-}(t)$, depending on $j$,
where
$$
\begin{array}{l}
p_{i,1,\pm}(t)=t^2+2(m_i\pm\alpha_{i,1})t\mp 2\alpha_{i,1}m_0-2m_im_0-m_0-1-2m_i\mp 2\alpha_{i,1},\\
p_{i,2,\pm}(t)=t^2+2(m_i\pm\alpha_{i,2})t+m_0+1,\\
p_{i,3,\pm}(t)=t^2+2(m_i\mp\alpha_{i,3})t-m_0-1-2m_i\pm2\alpha_{i,3}.
\end{array}
$$
These factorizations correspond, by the proposition \ref{2cov}, to
the $2$-torsion points in $E_i(\bQ)$ with $x$-coordinate equal to
$r_{i,1}=m_0m_i$, $r_{i,2}=-m_0-m_i-1$ and $r_{i,3}=0$.
\end{lemma}

By the previous lemma and the proposition \ref{2cov}, one can
construct Galois covers of $C_I$ with Galois group
$(\bZ/2\bZ)^2$, depending on the choice of the subset $J\subset I$
above and the choice of $j_1,j_2\in\{1,2,3\}$. The coverings can be described as the
projectivization (and desingularization) of the curve in $\bA^5$
given by
$$
C'\,:\,\{y_{1,+}^2\!\!= p_{1,j_1,+}(t)\,,\,y_{1,-}^2\!\!= p_{1,j_1,-}(t)\,,\,y_{2,+}^2\!\!= p_{2,j_2,+}(t)\,,\,y_{2,-}^2\!\!=  p_{2,j_2,-}(t)\,\},
$$
which is a curve of genus $17$, along with the map $\chi:C'\to
C_I$ defined as
$$\chi(t,y_{1,+},\,y_{1,-},\,y_{2,+},\,y_{2,-})=(t,y_{1,+}y_{1,-},\,y_{2,+}y_{2,-}).$$
These coverings can be defined over $\bQ$, although we choose to
show them in this form defined over the field
$\bQ(\alpha_{1,j},\alpha_{2,j})$, which is at most a biquadratic
extension of $\bQ$, in order to consider appropriate genus $1$
quotients of them.

Next, we choose one genus $1$ quotient of the form
$$H_{\pm,\pm}: z^2=p_{1,j_1,\pm}(t)p_{2,j_2,\pm}(t).$$
There are four such quotients, but depending on the degree of the
field $\bQ(\alpha_{1,j},\alpha_{2,j})$ there can be all of them
conjugates over $\bQ$, or to have two conjugacy classes if the
degree is 2, or all independent if the degree is 1.

For any element $\delta=(\delta_1,\delta_2)\in (\bQ^*)^2$, we
consider the twist $C'^{(\delta_1,\delta_2)}$ of the cover
$\chi$, given by
$$
C'^{(\delta_1,\delta_2)}\,:\,
\left\{
\begin{array}{lcl}
\delta_1y_{1,+}^2\!\!= p_{1,j_1,+}(t)\,&,&\,\delta_1y_{1,-}^2\!\!= p_{1,j_1,-}(t)\\
\delta_2 y_{2,+}^2\!\!= p_{2,j_2,+}(t)\,&,&\,\delta_2 y_{2,-}^2\!\!=  p_{2,j_2,-}(t)\,
\end{array}
\right\},
$$
along with the map
$$\chi^{(\delta_1,\delta_2)}(t,y_{1,+},\,y_{1,-},\,y_{2,+},\,y_{2,-})=(t,(y_{1,+}y_{1,-})/\delta_1^2,\,(y_{2,+}y_{2,-})/\delta_2^2).$$
We obtain 
$$C(\bQ) \subseteq \bigcup_{\delta \in \mathfrak{D}}
\chi^{(\delta)}(\{(t,y_{1,+},y_{1,-},y_{2,+},y_{2,-})\in
C'^{(\delta)}(\bQ(\alpha_{1,j_1},\alpha_{2,j_2}))\ : \ t\in \bP^1(\bQ) \}),$$
for some finite subset $\mathfrak{D} \subset (\bQ^*)^2$. The Proposition
\ref{2cov} allows us to describe the set $\mathfrak{D} $ in terms of
the Selmer groups of some isogenies. For any such
$\delta=(\delta_1,\delta_2)$, consider the quotients
$$H_{\pm,\pm}^{(\delta_1\delta_2)}: \delta_1\delta_2z^2=p_{1,j_1,\pm}(t)p_{2,j_2,\pm}(t)$$
which, in fact, only depend on the product $\delta_1\delta_2$. We
obtain
$$
\begin{array}{l}
\{t \in \bQ |\  \exists Y\in \bQ(\alpha_{1,j_1},\alpha_{2,j_2})^4 \mbox{ such
that } (t,Y)\in C'^{(\delta)}(\bQ(\alpha_{1,j_1},\alpha_{2,j_2}))\}\\
\qquad\qquad \subseteq  \{t \in \bQ |\  \exists w \in
\bQ(\alpha_{1,j_1},\alpha_{2,j_2}) \mbox{ such that }  (t,w) \in
H_{\pm,\pm}^{\delta}(\bQ(\alpha_{1,j_1},\alpha_{2,j_2}))\}.
\end{array}
$$
The following diagram illustrates, for a choice of $(j_1,j_2)$, all the curves
and morphisms involved in our problem:
$$
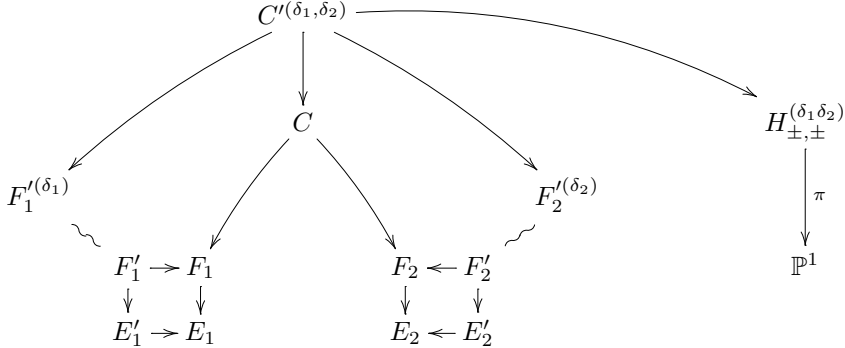
\xymatrix@R=0.7pc@C=0.8pc{
&             &            &      \ar@/_2mm/[dddlll] C'^{(\delta_1,\delta_2)} \ar@{->}[dd]    \ar@/^2mm/[dddrrr]   \ar@/^6mm/[ddrrrrrrr] &           &                           &      &       &     &   & \\
&                             &            &       &           &                           &      &    &        &   \\
&                     &            &   \ar@/_1mm/[ddl]     C      \ar@/^1mm/[ddr]                               &           &                           &   &  & & & H_{\pm,\pm}^{(\delta_1\delta_2)} \ar@{->}[dd]^\pi  \\
F_1'^{(\delta_1)}\ar@{~}[dr] &                                           &            &       &           &                          &   \ar@{~}[dl]  F_2'^{(\delta_2)}   &       &        & & \\
& F_1'  \ar@{->}[d] \ar@{->}[r]&  F_1 \ar@{->}[d]   &                                              &  F_2  \ar@{->}[d] &\ar@{->}[l] F_2' \ar@{->}[d] &  &  &  &&  \mathbb{P}^1   & \\
& E_1'  \ar@{->}[r]&  E_1    &                                              &  E_2   &\ar@{->}[l] E_2'  &      &   & & & &
}
$$
Note that the previous construction also depends on the choice of a subset $J=\{n_0,n_1,n_2\}\subset I$.

In the next lemma we describe a finite set $\mathfrak{S} \subset (\bQ^*)$
enough to cover all the possible values $t$ giving points of $C_I$, modulo the group of automorphisms $\Upsilon$.

\begin{lemma}\label{lema_delta} Let $I\subset\bN$ be a $5$-tuple. Fix a subset $J=\{n_0,n_1,n_2\}\subset I$ and $j_1,j_2\in\{1,2,3\}$.
For any $i=1,2$, denote by $\phi_i:E'_i\rightarrow E_i$ the $2$-isogeny
corresponding to the $2$-torsion point $(r_{i,j_i},0)\in E_i(\bQ)$, by $L=\bQ(\alpha_{1,j_1},\alpha_{2,j_2})$ and
by $\cS_L(\phi_i)$ a set of representatives in $L$ of the image
of $\Sel(\phi_i)$ in $L^*/(L^*)^2$ via the natural map. Finally,
denote by $\widetilde{\cS_L}(\phi_1)$ a set of representatives
of $\Sel(\phi_1)$ modulo the subgroup generated by the image of
the trivial points $\cT_I$ in this Selmer group. Consider the
subset $\mathfrak{S} \subset \bQ^*$ defined by
$$ \mathfrak{S} =\{\delta_1\delta_2 \ : \ \delta_1\in
\widetilde{\cS_L}(\phi_1), \delta_2\in \cS_L(\phi_2)\}.$$
Then, for any point $P=(t,y_1,y_2)\in C_I(\bQ)$, $\tau\in \Upsilon$ and $\delta\in \mathfrak{S}$ exist such that
$\tau(P)=(t',y_1',y_2')$ and $t'\in \bP^1(\bQ)$ such that $(t',w) \in H_{\pm,\pm}^{\delta}(L)$ for any sign $(\pm,\pm)$.
\end{lemma}

\begin{pf}
We have described at the previous paragraph that any point
$P=(t,y_1,y_2)\in C_I(\bQ)$ lift to a point in
$C'^{(\delta_1,\delta_2)}(L)$ for some $\delta_i \in
\Sel(\phi_i)$, for $i=1,2$. Hence determine a point in $
H_{\pm,\pm}^{\delta_1\delta_2}(L)$ with the first coordinate in
$\bP^1(\bQ)$.

If $P\in C_I(\bQ)$ has image $\delta_1$ in the Selmer group
$\Sel(\phi_1)$, then $\tau(P)$ has image $\delta_1\delta_{\tau}$,
if $\delta_{\tau}$ is the image of $\tau(T)$ in $\Sel(\phi_1)$ for
some trivial point $T$ with corresponding $\delta_1=1$. The previous sentence is true because the automorphisms belongs to $\Upsilon$ correspond to translation
by trivial points in the corresponding elliptic curve, if we fix (as we did) the zero point to be a trivial point (see Lemma 11 in
\cite{Xarles2012} for a proof in a special case). But the action of
$\Upsilon$ in $C_I(\bQ)$ is transitive on the set of trivial
points $\cT_I$.
\end{pf}

\

Now, the method allows us to conclude if we are able to compute,
for some choice of a subset $J=\{n_0,n_1,n_2\}\subset I$ and $j_1,j_2\in\{1,2,3\}$, and for any $\delta \in\mathfrak{S}$, all the points $(t,w) \in
H_{\pm,\pm}^{\delta}(\bQ(\alpha_{1,j_1},\alpha_{2,j_2}))$ with $t\in
\bQ$ for some choice of the signs $(\pm,\pm)$.

\

This last computation can be done in two steps as follows:
\begin{enumerate}
\item We first need to determine if there is some point in
$H_{\pm,\pm}^{\delta}(\bQ(\alpha_{1,j_1},\alpha_{2,j_2}))$. In the special
case $\delta=(1,1)$, the point at infinity always is a rational point. But, in
general, this curve (which will have points locally everywhere for
the $\delta$'s we choose) may have no rational points if it
represents an element of the Tate-Shafarevich group of its
Jacobian. We use the method described by Bruin and Stoll in \cite{Bruin-Stoll2009}. In particular, we have used their implementation in \verb|Magma| \cite{magma2.18-8} to determine if this happens.

\item Secondly, we will choose an isomorphism with its Jacobian
$\Jac(H_{\pm,\pm}^{\delta})$ and then we may use the elliptic curve Chabauty
technique as it was developed by Bruin at \cite{Bruin2003} to compute
this set if the rank of its group of
$\bQ(\alpha_{1,j_1},\alpha_{2,j_2})$-rational points is less than
the degree of $\bQ(\alpha_{1,j_1},\alpha_{2,j_2})$ over $\bQ$. We also need to determine a subgroup of finite index of this group to carry out the elliptic curve Chabauty method. Here, we have used the implementation in \verb|Magma| too.
\end{enumerate}

Hence we have $90$ possible choices of $J$, $j_1$ and $j_2$, and we
need to find one of them where we can carry out all these computations
for all the elements $\delta\in \mathfrak{S}$. In practice, we only consider the case where the field $\bQ(\alpha_{1,j_1},\alpha_{2,j_2})$ is
at most a quadratic extension of $\bQ$, essentially because of the
computation of the rank and/or a subgroup of finite index in
$\Jac(H_{\pm,\pm}^{\delta})(\bQ(\alpha_{1,j_1},\alpha_{2,j_2}))$ is
too expensive computationally for number fields of higher degree.

\subsection{The algorithm at work}
We have implemented in \verb|Magma V2.18-8| the algorithm developed above. In the following we describe this algorithm in a few examples. For these $5$-tuples $I\subset \bN$ we show how it works. In the case that the output of the algorithm is \verb|true| then we obtain $\cZ_I$, otherwise we give detailed information about the reasons why the algorithm does not work.

$\bullet$ $I=\{ 0, 1, 2, 4, 7\}$: this is the first case having no rang zero elliptic quotients. First, we need to choose a subset $J \subset I$, and two values $j_1,j_2\in\{1,2,3\}$ such that the field $L=\bQ(\alpha_{1,j_1},\alpha_{2,j_2})$ is of degree less or equal to $2$. The subset $J=\{1,4,7\}$ and the pair $(j_1,j_2)=(2,1)$ do the job. In this case $L=\bQ(\sqrt{10})$ and we have the following factorizations:
$$
\begin{array}{ll}
p_{1,2,+}(t)=t^2 - 10/3 t + 2 ,& p_{2,1,+}(t)=t^2 + 1/3(-2\sqrt{10} - 10) t + 1/3(4\sqrt{10} + 14),\\
p_{1,2,-}(t)= t^2 - 6t + 2 ,& p_{2,1,-}(t)=t^2 + 1/3(2\sqrt{10} - 10) t + 1/3(-4\sqrt{10} + 14).
\end{array}
$$
Note that in fact in this case we have $\bQ(\alpha_{1,2})=\bQ$. Next step is to compute the set $\mathfrak{S}$ (see Lemma \ref{lema_delta}). We have $\mathfrak{S}=\{1,2,3,6\}$. Now for any $\delta \in\mathfrak{S}$, we must compute all the points $(t,w) \in
H_{\pm,\pm}^{\delta}(\bQ(\sqrt{10}))$ with $t\in
\mathbb{P}^1(\bQ)$ for some choice of the signs $(\pm,\pm)$ where
$$
H_{\pm,\pm}^{\delta}\,:\,\delta w^2=p_{1,2,\pm}(t)p_{2,1,\pm}(t).
$$
For $\delta=1, 6$, we have that $\mbox{rank}_{\bZ}H^{\delta}_{+,+}(\bQ(\sqrt{10}))=1$ therefore we can apply elliptic curve Chabauty to obtain the possible values of $t$. For $\delta=1$ (resp. $\delta=6$) we obtain $t=\infty$ (resp. $t=0$). For the values $t=\infty$ and $t=0$ we obtain the trivial points $[1:\pm 1:\pm 1:\pm 1:\pm 1]\in C_I(\bQ)$. For $\delta=2, 3$, using the method described by Bruin and Stoll in \cite{Bruin-Stoll2009} we obtain $H^{\delta}_{-,+}(\bQ(\sqrt{10}))=\varnothing$.

The following table shows all the previous data, where at the last column appears the arithmetic progression attached to the corresponding $t$:
$$
\begin{array}{c}
\begin{array}{|c|c|c|c|c|c|}
\hline
\delta & \mbox{signs} & H^{\delta}_{\mbox{\tiny signs}}(L)=\varnothing? & \mbox{rank}_{\bZ}H^{\delta}_{\mbox{\tiny signs}}(L) & t & (q,a)\\
\hline
1 & (+,+) & \mbox{no} & 1 & \infty & (0,1)\\
\hline
2 & (-,+) & \mbox{yes} & - & - & -\\
\hline
3 & (-,+) & \mbox{yes} & - & - & -\\
\hline
6 & (+,+) & \mbox{no} & 1 & 0 & (0,1)\\
\hline
\end{array}\\[2mm]
I=\{0,1,2,4,7\}, J=\{1,4,7\}, (j_1,j_2)=(2,1), L=\bQ(\sqrt{10})
\end{array}
$$
Looking at the previous table, we obtain $C_I(\bQ)=\{[1:\pm 1:\pm 1:\pm 1:\pm 1]\}$; and therefore $z_I=0$ if $I=\{ 0, 1, 2, 4, 7\}$.

\

$\bullet$ $I=\{ 0, 1, 2, 5, 7\}$: this is the Rudin sequence. Let be $J=\{2,5,7\}$ and $(j_1,j_2)=(3,2)$. Then we have $L=\bQ(\sqrt{14})$ and $\mathfrak{S}=\{\pm 1,\pm 2,\pm 5,\pm 10\}$. The following table summarises all the computations made in this case:
$$
\begin{array}{c}
\begin{array}{|c|c|c|c|c|c|}
\hline
\delta & \mbox{signs} & H^{\delta}_{\mbox{\tiny signs}}(L)=\varnothing? & \mbox{rank}_{\bZ}H^{\delta}_{\mbox{\tiny signs}}(L)  & t & (q,a)\\
\hline
1 & (+,+) & \mbox{no} & 1 & \infty & (0,1)\\
\hline
-1 & (+,+) & \mbox{no} & 1 & - & -\\
\hline
2 & (+,-) & \mbox{no} & 1 & 3 & (24,1)\\
\hline
-2 & (+,-) & \mbox{no} & 1 &- & -\\
\hline
5 & (+,-) & \mbox{no} & 1 & 5/6 & (24,1)\\
\hline
-5 & (+,-) & \mbox{no} & 1 & - & -\\
\hline
10 & (+,+) & \mbox{no} & 1 & 0 & (0,1)\\
\hline
-10 & (+,+) & \mbox{no} & 1 & - & -\\
\hline
\end{array}\\[2mm]
{I=\{0,1,2,5,7\}, J=\{2,5,7\}, (j_1,j_2)=(3,2), L=\bQ(\sqrt{14})}
\end{array}
$$
We have that $C_I(\bQ)=\{[1:\pm 1:\pm 1:\pm 1:\pm 1], [1: \pm 5: \pm 7: \pm 11: \pm 13]\}$. That is, $\cZ_I=\{(24,1)\}$ for $I=\{ 0, 1, 2, 5, 7\}$.

\

$\bullet$ $I=\{0, 1, 3, 7, 8 \}$: this is an interesting example where there appear many values for $t$. Looking at the table below we obtain $\cZ_I=\{(120,1)\}$.
$$
\begin{array}{c}
\begin{array}{|c|c|c|c|c|c|}
\hline
\delta & \mbox{signs} & H^{\delta}_{\mbox{\tiny signs}}(L)=\varnothing? & \mbox{rank}_{\bZ}H^{\delta}_{\mbox{\tiny signs}}(L)  & t & (q,a)\\
\hline
1 & (+,+) & \mbox{no} & 1 & \begin{array}{c}\infty,1\\ 4,5/6\end{array} & \begin{array}{c}(0,1)\\(120,1)\end{array}\\
\hline
-1 & (+,+) & \mbox{no} & 1 & \begin{array}{c}0,3/2\\9/5,3/8\end{array} & \begin{array}{c}(0,1)\\(120,1)\end{array}\\
\hline
2 & (+,+) & \mbox{yes} & - & - & -\\
\hline
-2 & (+,+) & \mbox{yes} & - & - & -\\
\hline
\end{array}\\[2mm]
{I=\{0, 1, 3, 7, 8 \}, J=\{1,3,7\}, \{j_1,j_2\}=\{3,3\},L=\bQ(\sqrt{7})}
\end{array}
$$
In this case we have obtained $C_I(\bQ)=\{[1:\pm 1:\pm 1:\pm 1:\pm 1], [1: \pm 11: \pm 19: \pm 29: \pm 31]\}$.

\

$\bullet$ $I=\{0, 1, 4, 7, 8\}$: in this case we have at less two possible choices of $J$ and $(j_1,j_2)$ where the algorithm works obtaining $z_I=0$. In the first case $L=\bQ(\sqrt{2})$:
$$
\begin{array}{c}
\begin{array}{|c|c|c|c|c|c|}
\hline
\delta & \mbox{signs} & H^{\delta}_{\mbox{\tiny signs}}(L)=\varnothing? & \mbox{rank}_{\bZ}H^{\delta}_{\mbox{\tiny signs}}(L)  & t & (q,a)\\
\hline
1 & (+,+) & \mbox{no} & 1 & 0,\infty & (0,1)\\
\hline
3 & (+,+) & \mbox{no} & 1 & - & -\\
\hline
7 & (+,+) & \mbox{no} & 1 & - & -\\
\hline
21 & (-,+) & \mbox{no} & 1 &- & -\\
\hline
\end{array}\\[2mm]
{I=\{0, 1, 4, 7, 8\}, J=\{1, 4, 8\}, \{j_1,j_2\}=\{2,2\}, L=\bQ(\sqrt{2})}
\end{array}
$$
And what is more remarkable, also over $L=\bQ$ as the table below shows:
$$
\begin{array}{c}
\begin{array}{|c|c|c|c|c|c|}
\hline
\delta & \mbox{signs} & H^{\delta}_{\mbox{\tiny signs}}(L)=\varnothing? & \mbox{rank}_{\bZ}H^{\delta}_{\mbox{\tiny signs}}(L)  & t & (q,a)\\
\hline
1 & (+,-) & \mbox{no} & 0 & 1,\infty & (0,1)\\
\hline
2 & (+,+) & \mbox{yes} & - & - & -\\
\hline
-3 & (+,+) & \mbox{no} & 0 & 0,2 & (0,1)\\
\hline
-6 & (+,+) & \mbox{yes} & - & - & -\\
\hline
\end{array}\\[2mm]
{I=\{0, 1, 4, 7, 8\}, J=\{1, 4, 7\}, \{j_1,j_2\}=\{2,1\}, L=\bQ}
\end{array}
$$

\

$\bullet$ $I=\{0,3,5,6,10\}$: this is the second $5$-tuple where the algorithm does not work. The first one is $I=\{0,1,2,6,10\}$ and the reason is that 10 CPU hours was not enough to finish the computations for $I$.  In the following table appear all the subsets $J\subset I$ and pairs $(j_1,j_2)$ such that $L=\bQ(\sqrt{D})$ for some $D\in\bZ$. Note that in all the previous cases, $p_1(t)$ and $p_2(t)$ do not factorize over $\bQ$. Therefore is enough to check the signs $(+,+)$ and $(-,+)$. For $\delta=1$ we have computed an upper bound of the rank (denoted by $\mbox{rank}^*$) of the Mordell-Weil group of the Jacobians of the curves $H^{1}_{+,+}(L)$ and $H^{1}_{-,+}(L)$, where in all those cases is greater than $1$. Therefore we can not apply elliptic curve Chabauty and the algorithm outputs \verb|false|.
$$
\begin{array}{c}
\begin{array}{|c|c|c|c|c|}
\hline
J & \{j_1,j_2\}& D & \mbox{rank}^*_{\bZ}H^{1}_{+,+}(\bQ(\sqrt{D}))  & \mbox{rank}^*_{\bZ}H^{1}_{-,+}(\bQ(\sqrt{D}))\\
\hline
\{0,5,10\} &   \{2,3\} & -6 &  2 & 2 \\
\hline
\{0,3,6\} &   \{2,3\}& 10 & 2 & 2 \\
\hline
\{0,6,10\} &   \{2,3\}& -1 & 3 & 2 \\
\hline
\{0,3,5\} &   \{2,3\}& 2 & 2 & 3 \\[1mm]
\hline
\end{array}\\[2mm]
{I=\{0,3,5,6,10\}}
\end{array}
$$

$\bullet$ $I=\{ 0, 2, 4, 5, 11 \}$: this example shows one case where for all subsets $J\subset I$ of three elements and for all $j_1,j_2\in\{1,2,3\}$ we have that $L=\bQ(\alpha_{1,j_1},\alpha_{2,j_2})$ is a biquadratic extension of $\bQ$.

\

$\bullet$ Note that for all the $5$-tuples $I\subset \bN$ such that our algorithm has worked out we have obtained $z_I=0$ or $z_I=1$, except in the case $I=\{ 0, 13, 24, 33, 49 \}$. The table below shows that $\mathcal Z_I=\{(24,49),(-1,49)\}$, that is $z_I=2$.
$$
\begin{array}{c}
\begin{array}{|c|c|c|c|c|c|}
\hline
\delta & \mbox{signs} & H^{\delta}_{\mbox{\tiny signs}}(L)=\varnothing? & \mbox{rank}_{\bZ}H^{\delta}_{\mbox{\tiny signs}}(L)  & t & (q,a)\\
\hline
1 & (+,+) & \mbox{no} & 1 & \begin{array}{c}\infty,0\\ -12,-2/11\end{array} & \begin{array}{c}(0,1)\\(-1,49)\end{array}\\
\hline
6 & (+,+) & \mbox{no} & 1 & -12 & (-1,49)\\
\hline
10 & (+,-) & \mbox{no} & 0 & 2,12/11& (-1,49)\\
\hline
11 & (+,-) & \mbox{no} & 0 & 2,12/11& (-1,49)\\
\hline
14 & (+,-) & \mbox{no} & 1 & 12/11& (-1,49)\\
\hline
21 & (+,+) & \mbox{no} & 1 & \begin{array}{c}16/3\\ -12\end{array} & \begin{array}{c}(24,49)\\(-1,49)\end{array}\\
\hline
35 & (+,-) & \mbox{no} & 0 & 2,12/11& (-1,49)\\
\hline
154 & (+,-) & \mbox{no} & 0 & 2,12/11& (-1,49)\\
\hline
\end{array}\\[2mm]
{I=\{0, 13, 24, 33, 49  \}, J=\{0,13,24\}, \{j_1,j_2\}=\{2,2\},L=\bQ(\sqrt{165})}
\end{array}
$$
In this case we have obtained $C_I(\bQ)=\{[1:\pm 1:\pm 1:\pm 1:\pm 1], [49: \pm 36: \pm 25: \pm 16: \pm 0], [49: \pm 361: \pm 625: \pm 841: \pm 1225]\}$.

\section{Summary of the computations}\label{sec_comp_result}

One of the main objectives of this article is to prove Rudin's conjecture until $N=52$. For this purpose, we have developed a method based on the computation of the rational points of the curves $C_I$ attached to finite subsets $I\subset \bN$. Furthermore, at section \ref{preliminaries}, we have defined an equivalence relation on the finite subsets of $\bN$ such that for any pair of finite subsets $I,J$ such that $I\sim J$ we have that $C_I$ is $\bQ$-isomorphic to $C_J$, in particular $z_I= z_J$; and such that in any given equivalence class we have a primitive representant. Therefore, we restrict our computations to primitive subsets. Remember that we have attached to any finite subset $I\subset \bN$ a positive integer $n_I$, then we can introduce an ordering on equivalence classes of finite subsets of $\bN$.

First, let us consider the case of finite subsets $I\subset
\bN\cap\{0,\dots,51\}$ of cardinality $4$. There are $270725$ of
those subsets. But only $9077$ equivalence classes. We have proved
at section \ref{4squares} that the corresponding curves are elliptic
curve over $\bQ$.  The following table shows the number of curves
for a given rank:
$$
\begin{array}{|c||c|c|c|c|c|}
\hline
\mbox{rank} &  0 & 1 & 2 & 3 & 4 \\
\hline
\#\mbox{ curves} & 199 & 4692 & 3778 & 406 & 2\\
\hline
\end{array}
$$
If we restrict our attention to the symmetric case, we only have $402$ equivalence classes:
$$
\begin{array}{|c||c|c|c|}
\hline
\mbox{rank} &  0 & 1 & 2 \\
\hline
\#\mbox{ curves} & 190 & 191 & 2\\
\hline
\end{array}
$$

Next step it is to compute with subsets of $5$ elements. There are
$2598960$ of those subsets of $\{0,\dots,51\}$, $117449$ equivalence
classes. Then we remove all the subsets $I$ in the previous list with a subset $J$ such that $C_J(\bQ)$ is an elliptic curve of
rank $0$ and it only has $8$ torsion points, since in that case
$z_I=0$. After this sieve, $111338$ subsets remain. Now, at section
\ref{5squares}, given a subset $I\subset \bN$ with $5$ elements we
have developed a method that allows in some case to determine
$C_I(\bQ)$. The method consists, first to choose a subset
$J\subset I$ of three elements and $j_1,j_2\in\{1,2,3\}$. There are
$90$ possible choices. Secondly, to compute the finite set
$\mathfrak{S}$. Afterwards, compute for any $\delta
\in\mathfrak{S}$, all the points $(t,w) \in
H_{\pm,\pm}^{\delta}(\bQ(\alpha_{1,j_1},\alpha_{2,j_2}))$ with $t\in
\mathbb{P}^1(\bQ)$ for some choice of the signs $(\pm,\pm)$. This
method has worked out in $26589$ genus $5$ curves $C_I$. For those,
there are $26165$ cases such that $C_I(\bQ)=\mathcal{T}_I$ and $424$
cases such that $C_I(\bQ)\neq\mathcal{T}_I$. For the remaining
cases, $84749$, our method does not work for different reasons.
Notice that for a fixed $\delta\in\mathfrak{S}$ and a choice of
$(\pm,\pm)$ there are different reason making that our method will
not work.  First, we have bounded our computations for the case
where the field $\bQ(\alpha_{1,j_1},\alpha_{2,j_2})$ is at most a
quadratic extension of $\bQ$, since the algorithms on \verb|Magma|
that we are going to use are better implemented than in general
number fields. There are $34548$ cases where all the $90$ possible
choices give biquadratic fields. For the remaining cases, there are $1033$ such that \verb|Magma| crashed for some unknown reason or there had not been enough time (maximum of 10 CPU hours);  and there
are $49168$ cases where we know the reasons because our method has
not worked for them. For a given case, we need to decide if
$H_{\pm,\pm}^{\delta}(\bQ(\alpha_{1,j_1},\alpha_{2,j_2}))$ is empty
or not. Then the first reason such that our method does not work is:
\begin{itemize}
\item[(BS)] \verb|Magma| does not determine if $H_{\pm,\pm}^{\delta}(\bQ(\alpha_{1,j_1},\alpha_{2,j_2}))$ is empty or not.
\end{itemize}
Now assuming that we have computed a rational point on $H_{\pm,\pm}^{\delta}(\bQ(\alpha_{1,j_1},\alpha_{2,j_2}))$. Then two more reason may occur:
\begin{itemize}
\item[(Rank)] An upper bound of $\mbox{rank}_\bZ\,\jac{H_{\pm,\pm}^{\delta}}(\bQ(\alpha_{1,j_1},\alpha_{2,j_2}))$ is greater than one. Then, in principle, we can not use the elliptic curve Chabauty method.
\item[(noMW)] \verb|Magma| does not determine a subgroup of finite index on the elliptic curve $\jac{H_{\pm,\pm}^{\delta}}(\bQ(\alpha_{1,j_1},\alpha_{2,j_2}))$.
\end{itemize}
Notice that more than one reason could happen for a given $I$ making that any of the $90$ possible choices do not compute $C_I(\bQ)$ by our method. The next table shows the number of cases for the corresponding reasons:
$$
\begin{array}{lllllll}
\mbox{(Rank)}&: 37394 &  &\!\!\!\! \mbox{(Rank)+(noMW)}&:988 &\\
\mbox{(BS)}&: 630& & \!\!\!\! \mbox{(BS)+(Rank)}&:8526 &  \mbox{(Rank)+(noMW)+(BS)}\,:1523\\
\mbox{(noMW)}&: 11 & &\!\!\!\! \mbox{(noMW)+(BS)}&:96  &
\end{array}
$$
 
 \
 
All these computations ($110305$ $5$-tuples such that the algorithm has finished in less than $10$ CPU hours) took around 68 days of CPU time on a MacPro4.1 with 2 x 2.26 GHz Quad-Core Intel Xeon. 

\

The first case where we have not been able to determine $C_I(\bQ)$ is $I=\{0,1,2,6,10\}$ since $10$ CPU hours was not enough. The second one is $I=\{0,3,5,6,10\}$. In this case, our method does not work since for all the elliptic quotients defined over quadratic fields the upper bound for the
rank is greater than $1$.


\section{Consequences and comments}\label{sec_consequences}
The main goal of this article is to give new evidences for the
Rudin's Conjectures. First, given a positive integer $N\ge 6$ the
Strong version asserts that $Q(N)=Q(N;24,1)$. Our strategy to
prove this conjecture is recursively, that is, if we know $Q(N)$
for some $N$ then we attempt to compute $Q(N+1)$. We have that
$Q(N)\le Q(N+1)\le Q(N)+1$. Therefore we must compute $\cZ_I$ for
any $I\subset \{0,\dots,N\}$ such that $\#I=Q(N)+1$. Note that if
$z_I=0$ for any such tuples $I$, then $Q(N+1)=Q(N)$. Otherwise
$Q(N+1)=Q(N)+1$.

In section \ref{first_cases} we have proved $Q(6)=Q(7)=4$, and $Q(8)=5$
since $Q(8;24,1)=5$. Following the same strategy, that is with the
computations of subsets of $4$ integers, we even prove that
$Q(9)=Q(10)=Q(11)=5$. But it is not enough to show that $Q(12)=5$,
since for $I=\{0,1,2,5,9,11\}$ all the genus $1$ quotients
attached to subsets of $I$ of four elements have positive rank.
However, by using the methods in section \ref{5squares}, we prove
that for $J=\{ 0, 1, 2, 9, 11 \}$ we have $z_J=0$, therefore
$z_I=0$.

Now, in the general case, the strategy we followed was to consider
all the primitive subsets $I$ of 5 elements in $\{0,\dots,51\}$
where we are not able to compute $C_I(\bQ)$, either using the
genus 1 quotients or by the methods in section \ref{5squares}, as
we have described in the section \ref{sec_comp_result}. Using this list we
recursively compute the list ${\mathcal NC}(k)$ of all the
primitive subsets $I$ of $k$ elements, $k\ge 6$, such that we are
not able to compute $C_I(\bQ)$, by finding all the primitive
subsets $I$ of $k$ elements whose subsets of $k-1$ elements are
equivalent to a subset in ${\mathcal NC}(k-1)$ (See table
\ref{tab:notdetermined}). Note that we determined $C_I(\bQ)$ for
all the subsets of $\{0,\dots,51\}$ with more than $10$ elements.
\begin{table}[h]
\centering
\begin{tabular}{|c|l|c|}
\hline
$ k$ & $I \subset \{0,\dots,51\}$ & number of $I$\\
\hline
5 &  \{0, 1, 2, 6, 10\} & 84749 \\
\hline
6 & \{0, 1, 2, 7, 12, 15\}\footnotemark[5]   & 289752 \\
\hline
7 & \{0, 1, 6, 8, 11, 19, 23\} & 299855 \\
\hline
8 & \{0, 1, 3, 11, 17, 22, 23, 30\} & 69241 \\
\hline
9 & \{0, 2, 4, 13, 14, 19, 30, 33, 41\} & 2082 \\
\hline 
10 & \{0, 2, 7, 14, 17, 24, 37, 40, 43, 48\} & 2 \\
\hline
\end{tabular}
\caption{We list the first primitive subsets (in the natural order
explained in section \ref{preliminaries}) with $k$ elements that
we are not able to determine $C_I(\bQ)$, together with the number
of such subsets.} \label{tab:notdetermined}
\end{table}\footnotetext[5]{Note that $I=\{0,1,2,7,12,15\}\subset \mathcal{S}_{16}(24,1)$. Therefore  $(24,1)\in \mathcal{Z}_I$, but we are not able to compute the exact value of $z_I$.}

Furthermore, using now the subsets $I$ with $5$ elements where we
explicitly determine $C_I(\bQ)$ and containing other points apart
from the trivial ones, we explicitly compute, for $N\ge 8$, all the
arithmetic progressions $(q,a)$ such that $\#\cS_{N}(q,a)=Q(N)$
except\footnotemark[4] for $N=11,12$ \footnotetext[4]{For the
$5$-tuples $\{0,1,2,6,10\},\{0,3,5,6,10\},  \{0,2,4,5,11\},\{0,2,5,7,11\}, \{0,1,5,8,11\}$ and
$\{0,1,6,8,11\}$ we have not been able to compute the rational
points of the corresponding genus $5$ curve.}. In the table
\ref{tab:determined} we summarize these results.

\begin{table}[t]
\centering
\begin{tabular}{|c|c|c|c|l|}
\hline
$ k$ & $\mathcal{GP}_k $ & $N$ & $Q(N)$ & Arithmetic Progressions $(q,a)$\\
\hline
\multirow{4}{*}{$-2$}  & \multirow{4}{*}{$7$} &  $8$&\multirow{4}{*}{$5$} & $(24,1)$ \\
 &  &$9-10$ &   & $(24,1)$, $(120,1)$ \\
 &  & $11\footnotemark[4]$&  & $(24,1)$, $(120,1)$, $(8,1)$ \\
 &  & $12\footnotemark[4]$ &  & $(24,1)$,  $(120,1)$, $(8,1)$, $(24,25)$,\\
 &  &  &  &  $\qquad\qquad\quad\qquad(120,49)$, $(40,1)$,$(168,1)$\\
\hline
\multirow{2}{*}{$3$}  & \multirow{2}{*}{$12$} &$13-14$  & \multirow{2}{*}{$6$} & $(24,1)$ \\
                                 &                                   &      $15$  &      & $(24,1)$, $(24,25)$, $(120,1)$ \\
\hline
 \multirow{4}{*}{$-3$}  & \multirow{4}{*}{$15$} &  $16-18$&\multirow{4}{*}{$7$} & $(24,1)$ \\
 &  &$19-20$ &   & $(24,1)$, $(120,49)$ \\
 &  & $21$&  & $(24,1)$, $(120,49)$, $(120,1)$ \\
 &  & $22$ &  & $(24,1)$, $(120,49)$, $(120,1)$, $(24,25)$, $(8,1)$ \\
\hline
 \multirow{3}{*}{$4$}  & \multirow{3}{*}{$22$} & $23$  & \multirow{3}{*}{$8$} & $(24,1)$ \\
 &  & $24-25$ & & $(24,1)$, $(120,49)$ \\
 &  &  $26$ & & $(24,1)$, $(120,49)$, $(24,25)$ \\
\hline
 \multirow{3}{*}{$-4$}  &  \multirow{3}{*}{$26$} & $27-31$ & \multirow{3}{*}{$9$} & $(24,1)$ \\
 &  &  $32-34$ & & $(24,1)$, $(120,1)$ \\
 &  & $35$ & &  $(24,1)$, $(120,1)$, $(24,25)$ \\
\hline
 \multirow{2}{*}{$5$}  & \multirow{2}{*}{$35$} & $36-39$ & \multirow{2}{*}{$10$} & $(24,1)$ \\
 &  & $40$ & & $(24,1)$, $(24,25)$ \\
\hline
 \multirow{3}{*}{$-5$}  &  \multirow{3}{*}{$40$} & $41-49$ &  \multirow{3}{*}{$11$} & $(24,1)$ \\
 & & $50$ & & $(24,1)$, $(120,49)$ \\
 &  & $51$ & & $(24,1)$, $(120,49)$, $(24,25)$ \\
\hline
$6$ & $51$ & $52$ & $12$ & $(24,1)$ \\
\hline
\end{tabular}
\vspace{1mm}\caption{In the first column $k$ is an integer, in the second the
corresponding generalized pentagonal number $\mathcal{GP}_k $, in
the third an integer $N$, in the fourth $Q(N)$ and in the last one
appear the arithmetic progressions $qn+a$ with $\gcd(q,a)$
squarefree and $q> 0$ such that have $Q(N)$ squares for
$n\in\{0,\dots,N-1\}$.} \label{tab:determined}
\end{table}

The computations from the table \ref{tab:determined} allow us to
prove what we have called  Super--Strong Rudin's Conjecture up to
level $52$: let be $8 \le N=\mathcal{GP}_k+1\le 52$ for some
integer $k$, then $Q(N)=Q(N;q,a)$ with $\gcd(q,a)$ squarefree and
$q>0$ if and only if $(q,a)=(24,1)$.

\

We finish this section by discussing some points concerning the
number of non-constant arithmetic progressions having their
squares in a subset $I \subset \{0,\dots,N\}$ with $\# I\ge
5$. One consequence of our computations is that, for the subsets
$I$ of $\{0,\dots,52\}$ with $\# I\ge 5$ that we are able to
compute $z_I$, we have obtained that $z_I\le 1$, except for one case where $z_I\le 2$. But it is easy to see that
this is not true in general.

\begin{lemma}  Consider $a_i,q_i\in\bZ$. Then:
\begin{enumerate}
\item If $q_1q_2$ is not a square, then the set
$\cS(q_1,a_1^2)\cap \cS(q_2,a_2^2)$ is infinite.
\item  $\cS(q_1,a_1^2)\cap \cS(q_2,a_2^2)\cap \cS(q_3,a_3^2)$ is finite.
\item If the Bombieri-Lang conjecture is true, there exists an $r$ such that, for
any set of $r$ pairs $(q_i,a_i)$ of coprime integers, $\bigcap_{i=1}^r \cS(q_i,a_i^2)$ has at most 4 elements.
\end{enumerate}
\end{lemma}

\begin{pf} The set $\cS(q_1,a_1^2)\cap \cS(q_2,a_2^2)$ can be
described also by the set of integer solutions of the equation
$$x_1^2-q_1q_2x_2^2=q_2^2a_1^2-q_1q_2a_2^2,$$
which is a Pell type equation with a solution. Hence it has an
infinite number of solutions.

In the case we have three pairs, we look for integer solutions of an
equation giving a genus 1 curve, so it has a finite number of them
by Siegel's Theorem.

If we have more than three pairs, the resulting curve will be of genus
bigger than 1. So, suppose we have $r$ pairs such that
$J=\bigcap_{i=1}^r \cS(q_i,a_i^2)$ has more than 4 elements, so
there is a subset $I\subset J$ with 5 elements in it. This means
that the corresponding curve $C_I$ will have genus $5$, and with
$\#C_I(\bQ)\ge 16r+8$ (and, if $q_i\ne 0$ for all $i \in I$, in fact
$\ge 16(r+1)$). But thanks to the results from \cite{Caporaso-Harris-Mazur1997}, the
Bombieri-Lang conjecture implies there is an absolute bound for the
number of rational points of genus 5 curves over $\bQ$. Hence such
$r$ is bounded by above.
\end{pf}

\begin{exmp}\label{ej48} Using the ideas of the previous lemma, it is easy to construct one parametric
families of subsets $I \subset \bN$ with 5 elements along with two different and non-constant arithmetic progressions taking
squares in $I$. For example, for any integer $s>1$, we have that
$$\cS(s-1,1)\,\cap\, \cS(s+1,1)\supset \{0,4s, 4s(4s^2-1),8s(8s^4-6s^2+1), 8s(32s^6-40s^4+14s^2-1)\}.$$

As a consequence, we have built a one-parametric family of genus
$5$ non-hyperellip\-tic curves (of the form $C_I$) having at least
$3\cdot16=48$ points.
\end{exmp}

\begin{rem}\label{elkies} If the Bombieri-Lang conjecture is true then thanks to the results from \cite{Caporaso-Harris-Mazur1997} we have that there exists a
bound $B(g,\bQ)$ such that any curve of genus $g$ defined over
$\bQ$ satisfying $\#C(\bQ)\le B(g,\bQ)$. For the special case
$g=5$, Kulesz \cite{Kulesz1998} found a biparametric family of
hyperelliptic curves of genus $5$ with $24$ automorphisms over
$\bQ$ with at least $96$ points such that specializing he is able
to found a genus $5$ hyperelliptic curve $C$ defined over $\bQ$
such that $\#C(\bQ)=120$.  For the non-hyperelliptic case we have
that the curve $C_I$ attached to a $5$-tuple $I\subset \bN$ is of
genus $5$ and has $16$ automorphisms over $\bQ$. The example
\ref{ej48} shows a one-parametric family of genus $5$
non-hyperellip\-tic curves with at least $3\cdot16=48$ points.
Furthermore, we have found the following curves attached to
$5$-tuples $I\subset \bN$ such that $\#C_I(\bQ)\ge 5\cdot 16=80$. For
this search, we have looked for $5$-tuples such that have points
corresponding to $\mathcal{S}( 24b, a)$ with $a=1+24k$ square for
some $b,k\in\bN$. The table \ref{tab:bigexamples} shows the
results we have obtained.

\begin{table}[h]
\centering
\begin{tabular}{|c|cccc|}
\hline
$I$ & \multicolumn{4}{c|}{Arithmetic progression $(q,a)$ such that $I\subset\mathcal{S}(q,a)$}\\
\hline
$\{ 0, 2, 13, 23, 2233 \}$ & $( 240, 1369 )$ & $ ( 72, 25 )$ &  $( 120, 3481 )$ &  $( 168, 625 )$\\
\hline
$\{ 0, 5, 19, 70, 1020 \}$ & $( 72, 1 )$ &  $( 120, 2209 )$ &  $( 552, 961 )$ & $( 24, 169 )$\\
\hline
$\{ 0, 5, 33, 70, 1183 \}$ & $( 1344, 169 )$ &  $( 72, 1849 )$ &  $( 816, 961 )$ &  $( 24, 169 )$\\
\hline
$\{ 0, 17, 52, 147, 290 \}$ & $( 120, 1681 )$ & $( 96, 49 )$ &  $( 24, 961 )$ &  $( 264, 2401 )$\\
\hline
\end{tabular}
\vspace{1mm}
\caption{Some $5$-tuples $I$ with $z_I\ge 4$.}
\label{tab:bigexamples}
\end{table}
Note that the Bombieri-Lang conjecture implies that, for $k\ge 5$, a constant $c(k)$ should exist such that $z_I \le c(k)$ for
all $I\subset \bN$ with $\# I=k$.  In particular, $z_I \le c(5)$ for all $I\subset \bN$. The previous examples show that $c(k)\ge 2$ and $c(5)\ge 4$ (but we believe $c(5)> 4)$.
\end{rem}

\

{\bf Acknowledgements.}  We would like to thank Noan Elkies for useful discussion about the remark \ref{elkies} and to Nils Bruin for fixing a bug in the elliptic curve Chabauty \verb+magma+ routine that appeared on Magma v2.18-7.  

\

{\bf Data:} All the \verb+Magma+ and \verb+Sage+ sources are available on the first author's webpage.




\end{document}